\documentclass[11pt,reqno]{article}
\usepackage{amsmath} 
\usepackage{amssymb}
\usepackage{amsfonts}
\usepackage{blindtext}
\usepackage{graphicx}
\usepackage{xcolor}
\usepackage{caption}
\usepackage{subcaption}
\usepackage[leftcaption]{sidecap}

\usepackage{enumitem}

\usepackage{framed, amsthm, tikz, pgfplots}
\usepackage{anysize}

\usetikzlibrary{arrows}

\usetikzlibrary{decorations.markings}

\begin{document}

\def\nexto{\kern -0.54em}

\def\R{{\mathbb R}}
\def\T{{\mathbb T}}
\def\S{{\mathbb S}}
\def\C{{\mathbb C}}
\def\Z{{\mathbb Z}}
\def\N{{\mathbb N}}
\def\H{{\mathbb H}}
\def\B{{\mathbb B}}

\def\T{{\bf T}}
\def\diam{\mbox{\rm diam}}
\def\rr{{\cal R}}
\def\mt{{\Lambda}}
\def\e{\emptyset}
\def\sd{{\bf S}^{d-1}}
\def\stwo{{\bf S}^{2}}
\def\so{{\bf S}^1}
\def\dQ{\partial Q}
\def\dk{\partial K}
\def\endofproof{{\rule{6pt}{6pt}}}
\def\ts{\tilde{\sigma}}
\def\tr{\tilde{r}}
\def\be{\begin{equation}}
\def\ee{\end{equation}}
\def\beqn{\begin{eqnarray*}}
\def\eeqn{\end{eqnarray*}}

\def\ms{\medskip}
\def\beqn*{\begin{eqnarray*}}
\def\eeqn*{\end{eqnarray*}}
\def\endofproof{{\rule{6pt}{6pt}}}
\def\nv{\nabla \varphi}
\def\wuloc{W^u_{\mbox{\footnotesize\rm loc}}}
\def\wsloc{W^s_{\mbox{\footnotesize\rm loc}}}

\def\dist{\mbox{\rm dist}}
\def\diam{\mbox{\rm diam}}
\def\pr{\mbox{\rm pr}}
\def\supp{\mbox{\rm supp}}
\def\Arg{\mbox{\rm Arg}}
\def\In{\mbox{\rm Int}}
\def\Im{\mbox{\rm Im}}
\def\span{\mbox{\rm span}}
\def\spec{\mbox{\rm spec}\,}
\def\Re{\mbox{\rm Re}}
\def\var{\mbox{\rm var}}
\def\conf{\mbox{\footnotesize\rm const}}
\def\Conf{\mbox{\footnotesize\rm Const}}
\def\Lip{\mbox{\rm Lip}}
\def\con{\mbox{\rm const}\;}

\def\saa{\Sigma_A^+}
\def\sa{\Sigma_A}
\def\san{\Sigma^-_A}

\def\kmax{\kappa_{\max}}
\def\kmin{\kappa_{\min}}
\def\ecc{\mbox{\rm ecc}}
\def\tB{\widetilde{B}}
\def\hh{{\cal H}}
\def\thh{\widetilde{\cal H}}
\def\hE{\widehat{E}}
\def\ep{\epsilon}

\def\di{\displaystyle}
\def\pr{\mbox{\rm pr}}
\def\dK{\partial K}
\def\tpsi{\tilde{\psi}}
\def\tq{\tilde{q}}
\def\tsi{\tilde{\sigma}}
\def\ts{\tilde{s}}
\def\tW{\widetilde{W}}
\def\bs{\bigskip}
\def\tE{\widetilde{E}}
\def\tX{\widetilde{X}}
\def\beqn{\begin{eqnarray*}}
\def\eeqn{\end{eqnarray*}}
\def\yo{y^{(0)}}
\def\tY{\widetilde{Y}}
\def\tx{\tilde{x}}
\def\K{{\mathcal K}}
\def\T{{\mathcal T}}

\def\pp{{\mathcal P}}
\def\tgamma{\tilde{\gamma}}
\def\ot{(\omega,\theta)}
\def\overOm{\bar{\Omega}}
\def\spec{{\rm spec}}
\def\mm{{\mathcal M}}
\def\dk{\partial K}

\def\gen{{\rm generalized}\:}
\def\bic{{\rm bicharacteristic}}
\def\bics{{\rm bicharacteristics}}
\def\singsup{{\rm sing}\:{\rm supp}}
\def\Rnmin0{\R^n \setminus \{0\}}
\def\WF{WF}
\def\ff{{\mathcal F}}
\def\tu{\tilde{u}}
\def\tv{\tilde{v}}
\def\txi{\tilde{\xi}}
\def\teta{\tilde{\eta}}
\def\tSig{\widetilde{\Sigma}}
\def\tpsi{\tilde{\psi}}
\def\la{\langle}
\def\ra{\rangle}
\def\M{{\mathcal M}}
\def\aa{{\mathcal A}}
\def\tx{\tilde{x}}
\def\ty{\tilde{y}}
\def\Pm{P^{(m)}}

\def\tpi{\tilde{\pi}}
\def\ttheta{\tilde{\theta}}
\def\tOmega{\widetilde{\Omega}}
\def\Deone{\Delta^{(1)}}
\def\Detwo{\Delta^{(2)}}
\def\Dethree{\Delta^{(3)}}
\def\Defour{\Delta^{(4)}}
\def\tDeone{\widetilde{\Delta}^{(1)}}
\def\tDetwo{\widetilde{\Delta}^{(2)}}
\def\tDethree{\widetilde{\Delta}^{(3)}}
\def\tDefour{\widetilde{\Delta}^{(4)}}
\def\cc{{\mathcal C}}
\def\ta{\tilde{a}}
\def\tb{\tilde{b}}
\def\tM{\widetilde{M}}
\def\tf{\tilde{f}}
\def\Gr{\mbox{\rm Gr}}
\def\lan{\lambda^{(n)}}
\def\tn{t^{(n)}}
\def\En{E^{(n)}}
\def\Fn{F^{(n)}}
\def\lambdan{\lambda^{(n)}}
\def\tF{\widetilde{F}}
\def\tW{\widetilde{W}}
\def\ll{{\mathcal L}}
\def\tdelta{\tilde{\delta}}
\def\ti{\tilde{i}}

\noindent
  {\Large\bf Lyapunov Exponents for Open Billiards in the Exterior of Balls} 

\bigskip

 \noindent
{\bf Amal Al Dowais\footnote{Department of Mathematics, College of Science and Arts, Najran University, Najran, Saudi Arabia;
E-mail: amalduas@nu.edu.sa}}\\
{\footnotesize  Department of Mathematics and Statistics, University of Western Australia, Crawley 6009 WA, Australia\\
E-mail: amal.aldowais@research.uwa.edu.au}

\medskip

 \noindent
{\bf Luchezar Stoyanov}\\
{\footnotesize  Department of Mathematics and Statistics, University of Western Australia, Crawley 6009 WA, Australia\\
E-mail: luchezar.stoyanov@uwa.edu.au}

\bigskip

\noindent
{\bf Abstract.} In this paper we consider the billiard flow in the exterior of several (at least three) balls in $\R^3$
with centres lying on a plane. We assume that
the balls satisfy the no eclipse condition (H) and their radii are small compared to the
distances between their centres.
We prove that with respect to any Gibbs measure on the non-wandering set of the billiard flow
the two positive Lyapunov exponents are different: $\lambda_1 > \lambda_2 > 0$.

\medskip

\medskip 


\medskip


\section{Introduction}
\renewcommand{\theequation}{\arabic{section}.\arabic{equation}}

Let $K$ be a subset of ${\R}^{3}$  of the form 
$K = K_1 \cup \ldots \cup K_{k_0} ,$
where $k_0 \geq 3$ and all $K_i$ are balls with small radii $r_i > 0$ (not necessarily the same) and centres $O_i$ lying on a 
given plane $\Gamma$ in $\R^3$. Set 
$\Omega = \overline{{\R}^3 \setminus K}$. We assume that $K$ satisfies the following {\it no-eclipse condition}: 
$${\rm (H)} \quad \quad\qquad  
\begin{cases} \mbox{\rm for every pair $K_i$, $K_j$ of different connected components of $K$}\cr
\mbox{\rm the convex hull of $K_i\cup K_j$ has no common points with any other}\cr
\mbox{\rm  connected component of $K$. }
\end{cases}
$$

With this condition, the {\it billiard flow} $\phi_t$ defined on the {\it sphere bundle} 
$$S(\Omega) = \{ (q,v) : q \in \Omega\:, \: v \in \S^2 \}$$
in the standard way is called an  {\it open billiard flow}. 
It has singularities, however its restriction to the {\it non-wandering set} $\Lambda$ has only 
simple discontinuities at reflection points.  

\ms

{\bf Assumptions:} We will assume that for some constant $d_0 > 0$ we have  $d_0 < |O_iO_j|\leq 1$ for all $i \neq j$ and for each $i$
the radius $r_i$ of ball $K_i$ is so that $r_i \leq r_0$ for some constant $r_0 < d_0/3$. The condition (H) implies that
there exists a constant $c_1 > 0$ so that if $q_1,q_2,q_3$ are the successive reflection points of a billiard trajectory in 
$\Omega$ and $\nu(q_2)$ is the outward unit normal to $\dk$ at $q_2$, then for the angle $\varphi$ between 
$\nu(q_2)$ and $\frac{q_3-q_2}{\|q_3-q_2\|}$ we have $\cos \varphi \geq c_1$.  
Moreover, the configuration of our balls implies that there exists a global constant 
$c_2 < 1$ such that  if the reflection points $q_{1}$ and $q_{3}$  belong to
different components of $K$, then $\cos \varphi \leq c_2$. Thus, in the latter case always $0 < c_1 \leq \cos \varphi \leq c_2 < 1$.
In what follows we will assume that the constants $d_0$, $r_0$ and $c_2$ satisfy $c_2d_0 > 2 r_0$.

\bs

\noindent
{\bf Remark.} Clearly if $O_1, \ldots, O_{k_0}$ are arbitrary points in a plane $\Gamma$ so that no three of them lie on a line,
then choosing sufficiently small numbers $r_i > 0$, the balls $K_i$ with centres $O_i$ and radii $r_i$ ($i=1,\ldots, k_0$)
satisfy the condition (H) and the above assumptions.

\bs

Let $\mu$ be an arbitrary Gibbs measure (therefore an ergodic measure) for the shift map $\sigma$ on the related symbol
space $\Sigma_A$ with $k_0$ symbols $\{1,2, \ldots, k_0\}$  (see Sect. 2). Using the natural isomorphism between the non-wondering
set $M_0$ of the billiard ball map $B$ and $\Sigma_A$ which conjugates $B$ with $\sigma$ (see Sect. 2), we can regard $\mu$ as
a measure on $M_0$.
Let $\lambda_1 \geq \lambda_2   > 0$ be the positive Lyapunov exponents for the billiard ball map $B$ with respect to the measure $\mu$. 

\bs

In this note we prove the following

\bs

\noindent
{\bf Theorem 1.1.} {\it Under the above assumptions, the Lyapunov exponents $\lambda_1$ and $\lambda_2$ 
are distinct, i.e. $\lambda_1 > \lambda_2$}.

\bs

This result may be surprising to some extend. 
Chernov conjectured in \cite{Ch3} (see p.17 there; see also \cite{Ch2}) that
all positive Lyapunov exponents for a periodic Lorentz gas with spherical scatterers should be the same. We are not aware of
any positive or negative results concerning this conjecture. However, as Chernov mentioned in \cite{Ch3}, it had been shown both 
theoretically and numerically that for $3D$ random Lorentz gases with a random configuration of scatterers the two positive Lyapunov 
exponents are different (see \cite{LBS}, \cite{DP}  and \cite{BL}). There are probably some more recent results in this direction.

\section{Preliminaries}
\renewcommand{\theequation}{\arabic{section}.\arabic{equation}}

\subsection{Billiard ball map vs billiard flow}

For $x\in \Lambda$ and a sufficiently small $\epsilon > 0$ the (strong) {\it stable} and {\it unstable manifolds} 
of size $\epsilon$ for the billiard flow are defined by
$$\tW_\epsilon^s(x) = \{ y\in S(\Omega) : d (\phi_t(x),\phi_t(y)) \leq \epsilon \: \: \mbox{\rm for all }
\: t \geq 0 \; , \: d (\phi_t(x),\phi_t(y)) \to_{t\to \infty} 0\: \}\; ,$$
$$\tW_\epsilon^u(x) = \{ y\in S(\Omega) : d (\phi_t(x),\phi_t(y)) \leq \epsilon \:\: \mbox{\rm for all }
\: t \leq 0 \; , \: d (\phi_t(x),\phi_t(y)) \to_{t\to -\infty} 0\: \} .$$
The corresponding tangent unstable/stable bundles are defined by
$\tE^u(x) = T_x \tW_\epsilon^u(x)$ and $\tE^s(x) = T_x \tW_\epsilon^s(x)$. 
Set
$$M_1 = \{ (q,v) \in \dk \times \S^1 : \langle v, \nu_K(q) > 0 \} ,$$
where $\nu_K(q)$ is the {\it outward unit normal} to $\dk$ at $q$. If $x \in M_1$ is so that $\phi_t(x) \in M_1$
for some positive $t > 0$, we consider the smallest $t > 0$ with this property and define $B(x) = \phi_t(x)$.
The {\it non-wandering set} $M_0$ of the {\it billiard ball map} $B$ is the set of all $x\in M_1$ so that the billiard trajectory
$\{\phi_t(x)\}_{t\in \R}$ is bounded. On this set we have a well-defined homeomorphism $B : M_0 \longrightarrow M_0$ which is locally
the restriction of a smooth map $M_1 \longrightarrow M_1$.
For any $x = (q,v) \in M_0$ the stable/unstable manifolds for the billiard ball map $B$ are defined by
$$W_\epsilon^s(x) = \{ y\in M_1 : d (B^n(x),B^n(y)) \leq \epsilon \:\:  \mbox{\rm for all }
\: n \in \N \; , \: d (B^n(x),B^n(y)) \to_{n\to \infty} 0\: \}\; ,$$
$$W_\epsilon^u(x) = \{ y\in M_1 : d (B^{-n}(x), B^{-n}(y)) \leq \epsilon \:\: \mbox{\rm for all }
\: n \in \N \; , \: d (B^{-n}(x),B^{-n}(y)) \to_{n\to \infty} 0\: \} ,$$
where $\N$ is the set of positive integers.


Given $x = (q,v) \in M_0$  and a small $t_0 > 0$, set $y = (q+t_0 v, v)$. For a small $\ep > 0$ and a sufficiently small
$\delta > 0$ (depending on $\ep$) the map
\begin{equation}\label{eq:2.1}
\Psi :W^u_\delta (x) \longrightarrow \tW^u_\ep (y)
\end{equation}
such that $\Psi(z,w) = (z + t\, w, w)$ for all $(z,w) \in W^u (x)$, where $t = t(z,w) > 0$,
is a local diffeomorphism. In the same way there is a local diffeomorphism from $W^s_\delta (x)$ to $\tW^s_\ep (y)$. Moreover
$D\Psi (x) : T_xM_0 \longrightarrow T_y \Lambda$ is an isomorphism so that
$D\Psi(x) (E^u(x)) = \tE^u(y)$ and $D\Psi(x) (E^s(x)) = \tE^s(y)$.

It is known (see \cite{Si1}, \cite{Si2}) that $\tW^u_\ep(y)$ has the form $\tW^u_\ep(y) = \tY$, where 
$$\tY = \{ (p,\nu_Y(p)) : p\in Y\}$$
for some smooth $2$-dimensional surface  $Y$ in $\R^3$ containing the point $y$ such that $Y$
is strictly convex with respect to the unit normal field $\nu_Y$, i.e. the curvature of $Y$ is strictly positive.
In a similar way one can describe $\tW^s_\ep(y)$ using a strictly concave local surface.

Returning to the map (\ref{eq:2.1}), let $x_1 = (q_1,v_1) = B(x)$, where $q_1 = q + d \, v$ for some $d > 0$.
Define $y_1 = (q_1+ t_1 v_1, v_1)$ for some small $t_1 > 0$ (e.g. $t_1 = t_0$). Then there is a local diffeomorphism
(assuming again that $\ep > 0$ and $\delta  > 0$ are small)
$\Psi_1 :W^u_\delta (x_1) \longrightarrow \tW^u_\epsilon(y_1)$
defined as above. The unstable manifold $\tW^u_\ep(y_1) = \tY_1$ has the form  $\tY_1 = \{ (p_1,\nu_Y(p_1)) : p_1\in Y_1\}$
for some smooth $2$-dimensional surface $Y_1$ in $\R^3$ containing the point $y_1$ such that $Y_1$
is strictly convex with respect to the unit normal field $\nu_{Y_1}$. Setting $t = d+ t_1$, the following two diagrams
are commutative:

$$\def\normalbaselines{\baselineskip20pt\lineskip4pt \lineskiplimit3pt}
\def\mapright#1{\smash{\mathop{\longrightarrow}\limits^{#1}}}
\def\mapdown#1{\Big\downarrow\rlap{$\vcenter{\hbox{$\scriptstyle#1$}}$}}
\begin{matrix}
W^u_\delta(x) &\mapright{B}& W^u_{\delta}(x_1)\cr 
\mapdown{\Psi}& & \mapdown{\Psi_1}\cr \tW^u_\ep(y) = \tY&\mapright{\phi_t}& \tW^u_{\epsilon}(y_1) = \tY_1
\end{matrix}
 \qquad , \qquad
\def\normalbaselines{\baselineskip20pt\lineskip3pt \lineskiplimit3pt}
\def\mapright#1{\smash{\mathop{\longrightarrow}\limits^{#1}}}
\def\mapdown#1{\Big\downarrow\rlap{$\vcenter{\hbox{$\scriptstyle#1$}}$}}
\begin{matrix}
E^u(x) &\mapright{DB(x)}& E^u(x_1)\cr 
\mapdown{D \Psi}& & \mapdown{D \Psi_1}\cr E^u(y) &\mapright{D\phi_t(y)}& E^u(y_1)
\end{matrix} 
$$ 
with appropriate choices of the small constants $\delta > 0$ and $\ep > 0$.

\medskip

For an open billiard, the billiard ball map $B : M_0 \longrightarrow M_0$  is naturally isomorphic to a transitive subshift of finite type.
Here we briefly describe the natural conjugacy.

Let $k_0 \geq 3$ be as in Sect.1 and let $A = (A(i,j))_{i,j=1}^{k_0}$ be the $k_0\times k_0$ matrix 
such that $A(i,j) = 1$ if $i\neq j$ and $A(i,i) = 0$ for all $i,j$. Consider the  symbol space 
$$\sa = \{  (i_j)_{j=-\infty}^\infty : 1\leq i_j \leq k_0\; , \; A (i_j\; ,\; i_{j+1}) = 1 \:\: \mbox{ \rm for all } \: j\; \},$$
with the product topology and the {\it shift map} $\sigma : \sa \longrightarrow \sa$ 
given by $\sigma ( (i_j)) = ( (i'_j))$, where $i'_j = i_{j+1}$ for all $j$.
Given $0 < \theta < 1$, consider  the {\it metric} $d_\theta$ on $\sa$ defined by
$d_\theta(z,w) = 0$ if $z = w$ and $d_\theta(z, w) = \theta^m$ if $ z= (z_i)_{i\in \Z}$, $w = (w_i)_{i\in \Z}$,
$z_i = w_i$ for $|i| < m$ and $m$ is maximal with this property (see e.g. \cite{B} or \cite{PP} for general information on 
symbolic dynamics). Now define $R : M_0 \longrightarrow \sa$ by $R(x) = \{i_j\}_{j\in \Z}$ so that
$\pr_1(B^j(x)) \in K_{i_j}$ for all $j \in \Z$. It is well-known  (see e.g. \cite{PS}) that this is a well-defined bijection
which defines a conjugacy between $B : M_0 \longrightarrow M_0$ and $\sigma : \sa \longrightarrow \sa$, namely $R\circ B = \sigma \circ R$. 
Moreover, with an appropriate choice of $\theta \in (0,1)$, $R$ is a homeomorphism.

Let $\phi$ be a H\"older-continuous function on $\sa$ with respect to a metric $d_\theta$. The {\it Gibbs measure}
$\mu = \mu_\phi$ defined by $\phi$, is the unique $\sigma$-invariant invariant probability measure on $\sa$
for which there exist two constants $P$ (the {\it topological pressure} of $\phi$) and
$C>0$ such that, for any cylinder $C[i_0,\dotsc, i_{n-1}]$ and for any point $z$ in this cylinder we have
\begin{equation}
  C^{-1} \leq \frac{\mu (\cc[i_0,\dotsc, i_{n-1}])}{e^{S_n \phi(z) - nP}} \leq C,
\end{equation}
where $S_n\phi (z) = \phi(z) + \phi(\sigma(z)) + \ldots + \phi(\sigma^{n-1}(z))$ (see \cite{B} or \cite{PP} for details).
Such a measure $\mu$ is ergodic and there exist constants $C> 0$ and $\rho \in (0,1)$ such that for every
cylinder $\cc = \cc [i_0,\dotsc, i_{n-1}]$ of length $n$ we have $\mu(\cc) \leq C \, \rho^n$.

\subsection{Multiplicative Ergodic Theorem}

Let $f$ be an invertible  transitive subshift of finite type over a bilateral symbol space $X$ and let $\mu$ be an 
invariant probability measure. We will assume in addition that $\mu$ is ergodic.
Let $L$ be an invertible  continuous linear cocycle over $f$ acting on a continuous $\R^d$-bundle $E$ over $X$.  Thus, 
$L(x) : E(x) \longrightarrow E(f(x))$ is a linear map for each $x\in X$ and  
$$L^n(x) = L(f^{n-1}(x)) \circ L(f^{n-2}(x)) \circ \ldots \circ L(f(x)) \circ L(x) : E(x) \longrightarrow E(f^n(x)) $$
for every integer $n \geq 1$.

Given an integer $p = 1, \ldots,d$, let $\Gr_p(\R^d)$ be the {\it Grassman manifold}  of the linear subspaces  of $\R^d$ of dimension $p$
endowed with the usual distance $d(U,V)$ between subspaces $U,V \in \Gr_p(\R^d)$ defined by
$$d(U,V) =  \max \{ |\la u, w\ra | : u \in U, w\in V^\perp , \|u\| = \|w\| = 1\} . $$

The following is Oseledets' Multiplicative Ergodic Theorem stated under the above assumptions (see  e.g.  
\cite{Ar}, \cite{BP} or \cite{V}   for  related detailed exposition and proofs; see also \cite{BPS}).

\bs

\noindent
{\bf Theorem 2.1.}
(Multiplicative Ergodic Theorem). 
{\it There exists a subset $\ll$ of $X$ with  $\mu(\ll) = 1$ such that:}

\ms

(a) {\it  For all $x\in \ll$ there exists}
$\displaystyle N(x) = \lim_{n\to\infty} (L^n(x)^t L^n(x))^{1/2n} .$

\ms

(b) {\it There exist an integer $k \geq 1$  
such that  the operator (matrix) $N(x)$ has $k$ distinct eigenvalues $t_k(x) < \ldots < t_{1}(x)$ for all $x \in \ll$.  }

\ms

(c) {\it There exist numbers $t_k < t_{k-1} < \ldots < t_2 <  t_{1}$ such that  
$(\tn_i(x))^{1/n} \to t_i$ for all $x \in \ll$ and all $i = 1, \ldots, k$.}

\ms

(d) {\it For all $x\in \ll$ and every $j = 1, \ldots, k$ the dimension $\dim(\En_j(x)) = m_j$ is constant
and there exists $\lim_{n\to \infty} \En_j(x) = E_j(x)$ in $\Gr_{m_j}(\R^{d})$.} 

\bs

(e) {\it For all $x\in \ll$ and every $j = 1, \ldots, k$ we have
$$\lim_{n\to \pm\infty} \frac{1}{n} \log \|L^n(x) \cdot u\| = \log t_j \quad , \quad u \in E_j(x)\setminus \{0\} $$
for all $j = 1, \ldots, k$.}

\bs

The numbers $\log t_k <  \log t_{k-1} < \ldots  <  \log t_2 <  \log t_{1}$
are the   (distinct) {\it Lyapunov exponents} of the cocycle $L$ over $f$. 
The {\it Oseledets' subspaces} $E_i(x)$ are invariant with respect to $L$.

Since the billiard ball map $B : M_0 \longrightarrow M_0$ is conjugate to an invertible transitive subshift of finite type (see Sect. 2.1),
the above theorem applies to it. For open billiards in $\R^3$ either $k = 2$ or $k = 4$. When $k = 2$, the positive Lyapunov exponents are
$\lambda_1 = \lambda_2 = \log t_1$ and $E_1(x) = E^u(x)$ for all $x\in \ll$, while the negative Lyapunov exponents are $\lambda_3 = \lambda_4 = - \lambda_1$. 
When $k = 4$, there are two positive Lyapunov exponents $0 < \lambda_2 = \log t_2 < \lambda_1 = \log t_1$, and $E_1(x)$ and $E_2(x)$ are one-dimensional
subspaces of $E^u(x)$ for all $x \in \ll$, while the negative Lyapunov exponents are $\lambda_4 = - \lambda_1 < \lambda_3 = - \lambda_2 < 0$.

\subsection{Words' frequency in symbolic dynamics}

Let $\sa$ and $\sigma: \sa \longrightarrow \sa$ be as in Sect. 2.1. A sequence  $y = \{y_1, y_2, \ldots, y_p\}$
of symbols 
$$y_i \in G_0 = \{1, 2, \ldots, k_0\}$$ 
will be called {\it admissible} if $y_i \neq y_{i+1}$ for all $i = 1, \ldots, p-1$.
Let $w = w_0, w_1, \ldots, w_{k-1}$ ($k>1$) be a word (factor)
in $\sa$, i.e. an admissible  sequence of symbols $w_i \in G_0$. Assuming $p \geq k$, denote by
$$|(y_1,y_2, \ldots, y_p)|_w$$
the number of those $i \leq p-k$ so that $(y_i, y_{i+1}, \ldots, y_{i-k+1}) = (w_0, w_1, \ldots, w_{k-1})$.
That is, $|(y_1,y_2, \ldots, y_p)|_w$ is the number of times the word $w$ appears (one-to-one) as a part of $y$.
Given $z = (i_j)_{j\in \Z} \in \sa$, the {\it frequency} $f_w(z)$ of $w$ in $z$ is defined by
$$f_w(z) = \lim_{n\to\infty} \frac{|(i_{-n}, i_{-n+1}, \ldots,  i_{-1} , i_0 , i_1,  \ldots , i_{n-1} , i_n ) |_w}{2n+1} $$
whenever the limit exists. It is well-known (see e.g. Theorem 1.13 in \cite{Ber}) that for any given $w$ the frequency
$f_w(z)$ exists for almost all $z \in \sa$ with respect to ergodic invariant measures. In our particular situation this
gives the following. 

\bs

\noindent
{\bf Proposition 2.2.} {\it Let $\mu$ be an ergodic invariant measure on $\sa$ and let $w = w_0, w_1, \ldots, w_{k-1}$ ($k>1$) 
be an admissible  word. Let $\cc = \cc[w_0, w_1, \ldots, w_{k-1}]$ be the cylinder of length $k$ in $\sa$ defined by $w$.
Then there exists a subset $\Sigma''$ of $\sa$ with $\mu(\Sigma'') = 1$ such that
$f_w(z)$ exists for all $z \in \Sigma''$ and $f_w(z) = \mu(\cc)$ for all $z \in \Sigma''$.}

\bs

Stating the above slightly differently, it follows that under the assumptions in the Proposition for every $z = (i_j)_{j\in \Z}  \in \Sigma''$ 
there exists the limit
$$\tf_w(z) = \lim_{n\to\infty} \frac{|(i_0, i_1 , \ldots , i_{n-1} )|_w}{n} = \mu(\cc) ,$$
where $\cc = \cc[w_0, w_1, \ldots, w_{k-1}]$ as above.

Assume now that $\mu = \mu_{\phi}$ is a {\bf fixed Gibbs measure} on $\sa$ defined by a H\"older continuous function 
$\phi$ on $\sa$. In particular $\mu$ is ergodic (see e.g. \cite{B} or \cite{PP}). Then there exist constants $C > 0$ and $\rho\in (0,1)$ 
so that for every cylinder $\cc$ in $\sa$ of length $n$ we have $\mu(\cc) \leq C\, \rho^n$.

By Oseledets' Theorem, there exists a subset $\Sigma'$ of $\sa$ with $\mu(\Sigma') = 1$
so that for every $x \in M_0$ with $R(x) \in \Sigma'$ the formulae for
the Lyapunov exponents hold for vectors in $E^u_j(x)$:
\begin{equation}
\lim_{m\to \pm\infty} \frac{1}{m} \log \|dB^m(x) \cdot u\| = \lambda_j \quad , \quad u \in E_j(x)\setminus \{0\} .
\end{equation}
{\bf Fix a subset $\Sigma'$} of $\sa$ with this property.

Let $k = 2k'$ (for some $k' > 1$) be a fixed 
integer -- we will say later how large it should be chosen. Given any $p, q\in G_0$, $p \neq q$,
let $w_{p,q} = (p,q, p,q, \ldots, p,q)$ be the word of length $k = 2k'$ involving only the symbols $p$ and $q$. Being a cylinder
of length $k$, for $\cc_{p,q} = \cc[w_{p,q}]$ we have $\mu(\cc_{p,q}) \leq C \, \rho^{k}$. Set
$$\Omega_k = \bigcup_{p\neq q} \cc_{p,q} ,$$
where the union is taken over all (ordered) pairs of distinct elements of $G_0$. The number of these pairs is
$k_0(k_0-1)$, so taking $k > 1$ sufficiently large, we have
\begin{equation}
\tdelta = \tdelta(k) = \mu(\Omega_k) \leq C k_0 (k_0-1) \rho^k < \frac{1}{2} .
\end{equation}
We now {\bf fix an integer} $k = 2k'$ so that (2.4) holds.

By using the consequence of Proposition 2.2 mentioned after it, for every word $w_{p,q}$ as above there exists 
a subset $\Sigma''_{p,q}$ of $\Sigma'$ with $\mu(\Sigma''_{p,q}) = 1$ such that
\begin{equation}
\lim_{n\to\infty} \frac{|( i_0 , i_1 , \ldots , i_{n-1}) |_{w_{p,q}}}{n} = \mu(\cc_{p,q}) < \tdelta < \frac{1}{2}
\end{equation}
for every $z = (i_j)_{j\in \Z}  \in \Sigma''_{p,q}$.
Set
\begin{equation}
\Sigma_0 = \bigcap_{p\neq q} \Sigma''_{p,q} , 
\end{equation}
where the intersection is taken over all (ordered) pairs of distinct elements of $G_0$. Then $\Sigma_0 \subset \Sigma'$
and $\mu(\Sigma_0) = 1$.

\bs

\noindent
{\bf Corollary 2.3.} {\it For the fixed number $k$ with {\rm (2.4)} and the subset $\Sigma_0$ of $\Sigma'$ defined as above, we
have the following:
for every $z = (i_j)_{j\in \Z}   \in \Sigma_0$ there exists an integer $t_0 = t_0(z) \geq 0$ such that for every integer $t \geq t_0$ 
the number of those $j = 0,1, \ldots, (t-1) k$ so that
amongst the entries $i_{j}, i_{j+1}, \ldots, i_{j+k-1}$
there are at least three distinct numbers\footnote{I.e. this is not a sequence like $p,q,p,q,p,q,\ldots$, a repetiotion of two symbols only.}
is at least $t k/2 - k$. }

\bs

\noindent
{\it Proof.} Given $z = (i_j)_{j\in \Z}  \in \Sigma_0$, for every pair $p \neq q$ of numbers in $G_0$ there exists an integer
$t_{p,q} = t_{p,q}(\xi) \geq 0$ so that 
\begin{equation}
\frac{|( i_0 , i_1 , \ldots ,  i_{tk-1}) |_{w_{p,q}}}{tk}  < \frac{1}{2}
\end{equation}
for all $t \geq t_{p,q}$. Set 
$\di t_0 = \max_{p, q\in G_0\,, \, p\neq q} t_{p,q} .$
Consider an arbitrary  integer $t \geq t_0$. Let $t$ be the number of those $j = 0,1,\ldots, (t - 1)k$ so that
$(i_{j}, i_{j+1}, \ldots, i_{j+k-1}) = w_{p,q}$ for some $p\neq  q$ in $G_0$. Then (2.7) gives $t < t k/2$.
Thus, the number of the remaining numbers $j$ between $0$ and $(t - 1)k$ is $> tk/2-k$. For each of them
$(i_{j}, i_{j+1}, \ldots, i_{j+k-1}) \neq w_{p,q}$ for any $p\neq  q$ in $G_0$, therefore the sequence
$i_{j}, i_{j+1}, \ldots, i_{j+k-1}$ contains at least three different symbols from $G_0$.
This proves the corollary.
\endofproof

\section{Proof of Theorem 1.1}
\setcounter{equation}{0}

We will use the approach of Petkov and Vogel in \cite{PV} and the  formula they derived for a certain representation of the 
Poincar\'e map related to a periodic billiard trajectory. We follow Sect. 2.3 in \cite{PS}, where this formula is proved
in details. In  \cite{PV} (and Sect. 2.3 in \cite{PS}) the formula derived concerns a periodic trajectory, however a similar set-up
and similar calculations can be used for general billiard trajectories. This is what we do here.

One could also use the formula of Chernov and Sinai for the curvature operators on unstable manifolds using continuous
fractions (see e.g. \cite{Ch1}, \cite{Ch2}, \cite{Si2}; see also \cite{Ch3} and \cite{BCST}). However for the purposes of this 
paper we found the formula from \cite{PV} easier to use.

We follow the notation in Sect. 2.3 in \cite{PS} with some small changes.

\bs

Let $\mu = \mu_{\phi}$ be a Gibbs measure on $ \sa$ defined by a H\"older continuous function $\phi$ on $\sa$,
and let $\lambda_1 \geq \lambda_2 > 0$ be the Lyapunov exponents of the billiard ball map $B$ on the non-wandering set
$M_0$ with respect to $\mu$ (using here the representation map $R$ to identify $M_0$ with $\sa$, etc.; see Sect. 2.1).
As in Sect. 2.3, {\bf fix an integer $k \geq 1$ with (2.4) and define the subset $\Sigma_0$ of $\Sigma'$ by (2.6). }
 
{\bf Fix an arbitrary $\tx_0  = (p_0, v_0) \in M_0$} so that $R(\tx_0) = z \in \Sigma_0$. 
We may and will assume that $\tx_0$ is so that the billiard trajectory generated by it is never perpendicular to $\dk$ at
a reflection point.

Take a small $\ep > 0$ and consider the local unstable manifold $\tW^u_\ep(x_0)$ for $x_0 = (q_0, n_X(q_0))$, for some
convex local surface $X$ in the complement of $K$ with normal field $n_X$ so that 
$$q_0 = p_0 + t_0 n_X(q_0) \quad, \quad v_0 = n_X(q_0),$$
for some small $t_0 > 0$. That is $x_0$ is determined by the fixed point $\tx_0 \in M_0$.
Let $\gamma$ be the billiard trajectory in $\Omega$ determined by $x_0$ and let $q_1, q_2, \ldots$ be its 
successive reflection points.

Given $x = (q, n_X(q)) \in \tW^u_\ep(x_0)$,  set $q_0(x) = q$, $t_0(x) = 0$ and denote by $t_1(x), t_2(x), \ldots$ 
the times of the consecutive reflections of $\gamma(x)$ at $\dk$. Then 
$t_j(x) = d_0(x) + d_1(x) + \ldots + d_{j-1}(x)$, where $d_j(x) = \| q_{j+1}(x) - q_{j}(x)\|$, $0\leq j.$
Given $t \geq 0$, denote by $u_t(q_0)$ the {\it shift} of $q_0$ along the trajectory $\gamma(x)$ after time $t$, that is
$u_t(q_0) = \pr_1(\phi_t(x))$.
Set
$$X_t = \{ u_t(q) : q\in X\}\;.$$
When $u_t(q_0)$ is not a reflection point of $\gamma(x)$, then locally
near $u_t(q_0)$, $X_t$ is a smooth convex $2$-dimensional surface in $\R^3$ with "outward" 
unit normal given by the {\it direction} $v_t(q_0) = \pr_2(\phi_t(x))$ of $\gamma(x)$ at $u_t(q_0)$.

Fix for a moment $t > 0$ such that $t_m(x) < t < t_{m+1}(x)$ for some $m \geq 1$ (a large number; later $m \to \infty$),
and assume that $q_0(s)$, $0 \leq s\leq a$, is a $C^3$ curve on $X$ with $q_0(0) = q_0$ 
such that for every $s \in [0,a]$ we have $t_m(x(s) ) < t < t_{m+1}(x(s))$,
where $x(s) = (q_0(s), \nu_X(q_0(s)))$.  Assume also that $a > 0$ is so small that for all
 $j = 1,2,\ldots,m$ the reflection points $q_j(s) = q_j(x(s))$ belong to
the same boundary component $\dk_{i_j}$ for every $s\in [0,a]$.

We will now estimate $\| d\phi_t(x)\cdot \dot{q}_0(0)\|$, where $\dot{q}_0(0) \in T_{q_0}X$.

Denote by $\Pi_0$ the hyperplane in $\R^3$ passing through the point $q$ and orthogonal to $n_X(q)$.
We will actually identify $\Pi_0$ with $T_qX$ with the $0$ at the point $q$, so this will be a $2$-dimensional 
vector subspace of $\R^3$. Let $q(u)$, $u \in U_0$, be a smooth ($C^3$) parametrisation of $X$ in an open neighbourhood
of $q$ in $X$, where $U_0$ is an open disk in $\R^2$ with centre $0$ and such that $q(0) = q$.  Assuming $U_0$ is
sufficiently small, we can now use it to parametrise $\tW^u_\ep(x_0)$ locally near $x_0$. Namely we identify every 
$(u,v)$ in $\Pi_0\times \Pi_0$ close to $(0,0)$ (e.g. assume $(u,v) \in U_0\times U_0$) with the point
$(x(u) , w) \in \tW^u_\ep(x_0)$ so that $w = n_X(q) + v$. Set $\omega_0 = n_X(q)$.

Next, for any $i \geq 1$, denote by $\Pi_i$ the plane in $\R^3$ passing through the point $q_i$ and orthogonal to
the line $q_iq_{i+1}$, and by $\omega_i$ the unit vector determined by the vector $\overrightarrow{q_i q_{i+1}}$.
We assume that for each $i$ the hyperplane $\Pi_i$ is endowed with a linear basis such that $q_i = 0$.

Given any $ i \geq 0$ and a pair $(u,v) \in \Pi_i \times \Pi_i$ sufficiently close to $(0,0)$, let $\ell(u,v)$ be the
oriented line passing through $u$ and having direction $\omega_i + v$. Here we identify the point $v$
with the vector $v$. If $(u,v)$ is sufficiently close to $(0,0)$, then $\ell(u,v)$ intersects
transversally $\dk$ at some point $p = p(u,v)$ close to $q_{i+1}$. Let $\ell'(u,v)$ be the oriented line
symmetric to $\ell(u,v)$ with respect to the tangent plane $T_p(\dk)$ to $\dk$ at $p$, and let $u'$ be the
intersection point of $\ell'(u,v)$ with $\Pi_{i+1}$ (clearly such a point exists for $(u,v)$ close
to $(0,0)$). There is a unique $v'\in \Pi_{i+1}$ such that $\omega_{i+1} + v'$ has the direction
of the line $\ell'(u,v)$ (see Figure 1). Thus, we obtain a map
$$\Phi_{i+1} : \Pi_i \times \Pi_i \ni (u,v) \mapsto (u',v') \in \Pi_{i+1}\times \Pi_{i+1} ,$$
defined for $(u,v)$ in a small neighbourhood of $(0,0)$. The smoothness of this map follows from
the smoothness of the billiard flow. 

Just like we did with $\Pi_0$, let $\Pi_m$ be the plane in $\R^3$ passing through the point $u_t(q)$ and orthogonal to 
$n_{X_t}(u_t(q))$. Let $q_m(u)$, $u \in U_m$, be a smooth ($C^3$) parametrisation of $X_t$ in an open neighbourhood
of $u_t(q)$ in $X_t$, where $U_m$ is an open disk in $\R^2$ with centre $0$ and such that $q_m(0) = u_t(q)$.  
Again we will assume that $U_m$ is sufficiently small and then we can use it to parametrise $\tW^u_\ep(x(t))$ locally near 
$x(t)$, simply by identifying every  $(u_m,v_m)$ in $\Pi_m\times \Pi_m$ close to $(0,0)$ with the point
$(x(u_m) , w_m) \in \tW^u_\ep(x(t))$ so that $w_m = n_{X_t}(u_t(q)) + v_m$. Set $\omega_{m+1} = n_{X_t}(u_t(q))$.



\begin{center}
\begin{tikzpicture}[scale=1]

\tikzset{ kb1/.style={postaction={decorate, decoration={markings,mark=at position .5 with {\arrow{>};}}}}}


\node[circle, minimum size= 4cm] (A) at  (6.1, 0.3) { };

\draw (4,0.62) [rotate=40] arc[x radius=3, y radius=3, start angle=120, end angle=30] coordinate (b); 

\draw (-3,1.8) arc[x radius=3, y radius=3, start angle=25, end angle=130]; 

\node at (7.5, 2.2) {$\mbox{\rm \footnotesize $\partial{K}$}$};

\node at (-6.3, 3.2) {$\mbox{\rm \footnotesize $\partial {K}$}$};


\draw [blue, thick, postaction={decorate,  decoration={markings, mark=at position 0.13 with {\arrow{stealth}}}  }]   
(-4,3)--  (5,2) coordinate[pos=0.12] (a); 

\node at (-3, 3.15) {$\mbox{\rm \footnotesize $\omega_i$}$};

\draw [blue, thick, postaction={decorate,  decoration={markings, mark=at position 0.32
with {\arrow{stealth}}}  }]  (5, 2) -- (3.8,6) coordinate [pos=0.31] (d);

\node at (5.2, 3) {$\mbox{\rm \footnotesize $\omega_{i+1}$}$};


\draw  [green!70!black, thick, postaction={decorate,  decoration={markings, mark=at position 0.635 with {\arrow{stealth}}}  }](10, 4) -- (0.5, 0.2);

\draw [green!70!black, thick, postaction={decorate,  decoration={markings, mark=at position 0.71 with {\arrow{stealth}}}  }] (-39/9, 0) -- (-33/9, 6) coordinate[pos=0.7] (c);

\node at (-4, 5.6) {$\mbox{\rm \footnotesize $\Pi_i$}$};

\node at (7.5, 3.4) {$\mbox{\rm \footnotesize $\Pi_{i+1}$}$};


\draw [green!70!black, thick] (5,2) -- (4,3); 
\draw [green!70!black, thick] (-3.0,2.5) -- (-4, 3);
\draw [green!70!black, thick, stealth-]  (-3,2.5) -- (a) node [below, black] {$\mbox{\rm \footnotesize \, \, $v$}$};
\draw [green!70!black, thick, -stealth]  (d) node [left, black] { $\mbox{\rm \footnotesize$v'$\,}$} --(4,3) ; 


\draw [brown, thick] (A.150) -- (c) node [midway, above, black]{$\mbox{\rm \footnotesize$l(u,v)$}$};

\draw [brown, thick] (A.150) --(0,5) node [right, black] {$\mbox{\rm \footnotesize $l'(u,v)$}$};


\node at (-4.3, 3.8) {$\mbox{\rm \footnotesize $u$}$};

\node at (4.1, 1.85) {$\mbox{\rm \footnotesize $u'$}$};

\node[circle, fill, inner sep=0.7pt] at (5,2) { };
\node at (5.4, 1.8) {$\mbox{\rm \footnotesize $q_{i+1}$}$};

\node[circle, fill, inner sep=0.7pt, black] at (-4,3) { };
\node at (-4.2, 2.9) {$\mbox{\rm \footnotesize $q_{i}$}$};
\node[circle, fill, inner sep=0.7pt, brown] at (A.150) { };
\node at (5, 1.3) {$\mbox{\rm \footnotesize $p(u,v)$}$};

\end{tikzpicture}

\vspace{-1cm}
Figure 1: The map $\Phi_i$

\end{center}


\bs

Consider the (locally defined) composition 
$$\Phi = \Phi_m \circ \ldots \circ \Phi_1 : \Pi_0 \times \Pi_0  \longrightarrow \Pi_m \times \Pi_m ,$$
and the linear map\footnote{In the case of a periodic trajectory $\gamma$, $\Pm_\gamma$ is called the {\it Poincar\'e map} of 
$\gamma$. Here however we do not assume $\gamma$ to be periodic.}
$$\Pm_\gamma = d \Phi (0,0) :  \Pi_0 \times \Pi_0  \longrightarrow \Pi_m \times \Pi_m .$$

It follows from the above definitions that there exists a natural (local) conjugacy\footnote{The top map $\Phi$ is only defined locally, 
essentialy from $U_0\times U_0$ to $U_m\times U_m$.}
$$\def\normalbaselines{\baselineskip20pt\lineskip4pt \lineskiplimit3pt}
\def\mapright#1{\smash{\mathop{\longrightarrow}\limits^{#1}}}
\def\mapdown#1{\Big\downarrow\rlap{$\vcenter{\hbox{$\scriptstyle#1$}}$}}
\begin{matrix}
\Pi_0\times \Pi_0 &\mapright{\Phi}& \Pi_m \times \Pi_m\cr 
\mapdown{\Omega_0}& & \mapdown{\Omega_m}\cr \tW^u_\ep(x_0) = \tX &\mapright{\phi_t}& \tW^u_{\epsilon}(x_t(q)) = \tX_t
\end{matrix}
$$
which defines a conjugacy between the corresponding linearizations:
$$
\def\normalbaselines{\baselineskip20pt\lineskip3pt \lineskiplimit3pt}
\def\mapright#1{\smash{\mathop{\longrightarrow}\limits^{#1}}}
\def\mapdown#1{\Big\downarrow\rlap{$\vcenter{\hbox{$\scriptstyle#1$}}$}}
\begin{matrix}
\Pi_0 \times \Pi_0 &\mapright{\Pm_\gamma}& \Pi_m \times \Pi_m\cr 
\mapdown{D \Omega_0}& & \mapdown{D \Omega_m}\cr E^u(\tX) &\mapright{D\phi_t(x_0)}& E^u (\tX_t)
\end{matrix} 
$$ 
Moreover the smooth maps $\Omega_0$ and $\Omega_m$ have uniformly bounded derivatives, so 
there exists a constant $T> 0$ independent of $x_0$, $m$ and $t$ so that
$$\frac{1}{C} \, \|\Pm_\gamma (u,v)\| \leq \|D\phi_t(x_0)\cdot (D\Omega_0(0,0)\cdot (u,v))\| \leq C\, \|\Pm_\gamma (u,v)\| $$
for all $(u,v) \in \Pi_0\times \Pi_0$.
In particular when $u = 0$ we can naturally identify $D\Omega_0(0,v)$ with $v \in E^u(q_0)$, so we get
\begin{equation}
\frac{1}{C} \, \|\pr_2(\Pm_\gamma (0,v))\| \leq \|D\phi_t(x_0)\cdot v)\| \leq C\, \|\pr_2(\Pm_\gamma (0,v))\| 
\end{equation}
for all $v \in E^u(q_0)$ (with sufficiently small $\|v\|$), where we denote $\pr_2(u',v') = v'$ for all $(u',v') \in \Pi_i \times \Pi_i$
and all $i$.
Therefore, for the {\it Lyapunov exponent} $\lambda(v)$
at the point $x_0$ in the direction of $v$ we get
\begin{equation}
\lambda(v) = \lim_{m\to\infty} \frac{1}{m} \log \|\pr_2(\Pm_\gamma (0,v) ) \| .
\end{equation}

Next, we describe the  representation of $\Pm_\gamma$ given in Sect. 2.3 in \cite{PS}.

Set $d_i = \|q_{i-1} - q_i\|$ and let $\alpha_i = T_{q_i}(\dk)$ be the {\it tangent hyperplane} to $\dk$ at $q_i$.
Let $\sigma_i : \R^3 \longrightarrow \R^3$ be the {\it symmetry} with respect to $\alpha_i$, and let 
$\Pi'_i$ be the hyperplane passing through $q_i$ and orthogonal to $\omega_{i-1}$ ($i\geq 1$). Then $\Pi'_i$ is 
parallel to $\Pi_{i-1}$ and
$$\sigma_i(\omega_{i-1}) = \omega_i \quad, \quad \sigma_i(\Pi'_i) = \Pi_i $$
(see Figure 2).
We identify $\Pi'_i$ with $\Pi_{i-1}$ using the translation along the line $q_{i-1}q_i$ ($i \geq 1$).
Let 
$$\pi_i : \Pi'_i \longrightarrow \alpha_i$$
be the {\it projection} along the vector $\omega_{i-1}$.
For the outward {\it unit normal field} $\nu_i = \nu_i(q)$ to $\dk$ for $q\in \dk$, we have
$\theta_i = 2 \la \nu_i , \omega_i \ra > 0$. 
From the Assumptions in Sect. 1, there exists a constant $c_1 > 0$ independent of $q$, $x_0$ and $i$ so that 
$$\la \nu_i , \omega_i \ra \geq c_1 > 0$$
for all $i$. Thus, $\theta_i \geq 2c_1 > 0$ for all $i$.

The {\it second fundamental form} of $\dk$ at $q_i$ is defined by $\la G_i(\xi), \eta\ra $ for $\xi,\eta \in \alpha_i$, where 
$$G_i = d\nu_i(q_i) : \alpha_i \longrightarrow \alpha_i .$$
Since this is a symmetric form, there exists a unique symmetric linear map
\be 
\tpsi_i : \Pi_i \longrightarrow \Pi_i
\ee
such that
\be 
\la \tpsi_i \sigma_i (\xi') , \sigma_i (\eta') \ra = - 2 \la \omega_{i-1} , \nu_i \ra\,\la G_i(\pi_i(\xi')) , \pi_i(\eta') \ra 
= \theta_i \,\la G_i(\pi_i(\xi')) , \pi_i(\eta') \ra 
\ee
for all $\xi', \eta' \in \Pi'_i$; then $\xi = \sigma_i(\xi'), \,\eta = \sigma_i(\eta') \in \Pi_i$.

Since we identify the hyperplanes $\Pi_{i-1}$ and $\Pi'_i$ (using the projection along the line $q_{i-1}q_i$) we have
$\sigma_i(\Pi_{i-1}) = \Pi_i$. Consider the compositions
\be 
s_i = \sigma_i \circ \sigma_{i-1} \circ \ldots \circ \sigma_1 .
\ee
Then $s_i(\Pi_0) = \Pi_i$ and $s_i^{-1}(\Pi_i) = \Pi_0$. Define the symmetric linear map $\psi_i$ by
\be 
\psi_i = s_i^{-1}\circ  \tpsi_i \circ s_i : \Pi_0 \longrightarrow \Pi_0 .
\ee
Then we can write
$$
\left(
\begin{array}{cc}
I & d_i I\\
\psi_i  & I + d_i \psi_i
\end{array}
\right)
= \left(
\begin{array}{cc}
s_i^{-1} & 0\\
0  & s_i^{-1}
\end{array}
\right)
\left(
\begin{array}{cc}
I & d_i I\\
\tpsi_i  & I + d_i \tpsi_i
\end{array}
\right)
 \left(
\begin{array}{cc}
s_i & 0\\
0  & s_i
\end{array}
\right) ,
$$
so this defines a linear map $\Pi_0\times \Pi_0 \longrightarrow \Pi_0 \times \Pi_0$. 
Now the argument in Sect. 2.3 in \cite{PS} yields the following (see Theorem 2.3.1 in \cite{PS}):

\bs

\noindent
{\bf Lemma 3.1.} {\it Under the assumptions and conventions above, the map
$$\Pm_\gamma : \Pi_0 \times \Pi_0 \longrightarrow \Pi_m \times \Pi_m$$
is a linear symplectic map which has the following matrix representation
\be 
\Pm_\gamma =
 \left(
\begin{array}{cc}
s_m & 0\\
0  & s_m
\end{array}
\right) \cdot
\left(
\begin{array}{cc}
I & d_m I\\
\psi_m  & I + d_m \psi_m
\end{array}
\right)
\cdots
\left(
\begin{array}{cc}
I & d_1 I\\
\psi_1  & I + d_1 \psi_1
\end{array}
\right) ,
\ee
where $\psi_i$  is given by {\rm (3.6)} for all $i = 1,2 \ldots, m$.}

\bs

\begin{center}


\begin{tikzpicture}[scale=1.3]

\tikzset{ kb1/.style={postaction={decorate, decoration={markings,mark=at position .5 with {\arrow{>};}}}}}


\draw (3,0) arc[x radius=3.45, y radius=3.45, start angle=25, end angle=155];
\draw[red, thick, postaction={decorate,  decoration={markings, mark=at position 0.92 with {\arrow{stealth}}}  }]  (3,2)--(-3.2,2);
\draw[red, thick, -stealth ] (0,2)--(0,4);
\draw [blue, thick, postaction={decorate,  decoration={markings, mark=at position 0.9 with {\arrow{stealth}}}  }] (0, 2) --(3, 3);
\draw [blue, thick, postaction={decorate,  decoration={markings, mark=at position 0.3 with {\arrow{stealth}}}  }] (-3, 3) --(0, 2);
\draw [black, dashed, thick] (-0.3, 1.3) --(-2.8, 2);
\draw [black, dashed, thick] (-0.3, 2.8) --(-2.8, 2);
\draw [green!70!black, thick, postaction={decorate,  decoration={markings, mark=at position 0.63 with {\arrow{stealth}}}  }] (1, 5) -- (-1, -1);
\draw [green!70!black, thick, postaction={decorate,  decoration={markings, mark=at position 0.64 with {\arrow{stealth}}}  }] (1, -1) -- (-1, 5);
\node at (3, 2.7) {$\mbox{\rm \footnotesize $\omega_i$}$};
\draw[black,thick] (1.85,2.62) circle (0.02cm);

\node at (-2.4, 2.6) {$\mbox{\rm \footnotesize $\omega_{i-1}$}$};
\draw[black,thick] (-2.85,2.95) circle (0.02cm);

\node at (2.8, 1) {$\mbox{\rm \footnotesize $\partial{K}$}$};
\node at (-1.5, -0.5) {$\mbox{\rm \footnotesize $\Pi_{i-1}= \Pi'_i$}$};
\node at (1.1, -0.5) {$\mbox{\rm \footnotesize $\Pi_i$}$};
\node at (-0.1, 1.2) {$\mbox{\rm \footnotesize $\xi'$}$};
\node at (0, 4.1) {$\mbox{\rm \footnotesize $\nu_i$}$};
\node at (3.3, 2) {$\mbox{\rm \footnotesize $\alpha_i$}$};
\node at (-0.9, 2.9) {$\mbox{\rm \footnotesize $\xi= \sigma_i(\xi')$}$};
\node at (-3, 1.7) {$\mbox{\rm \footnotesize $\pi_i(\xi') =\tpi_i(\xi)$}$};
\node[circle, fill, inner sep=0.9pt, blue] at (0,2) {};
\node at (0.2, 1.84) {$\mbox{\rm \footnotesize $q_i$}$};

\end{tikzpicture}

\end{center}

\begin{center}
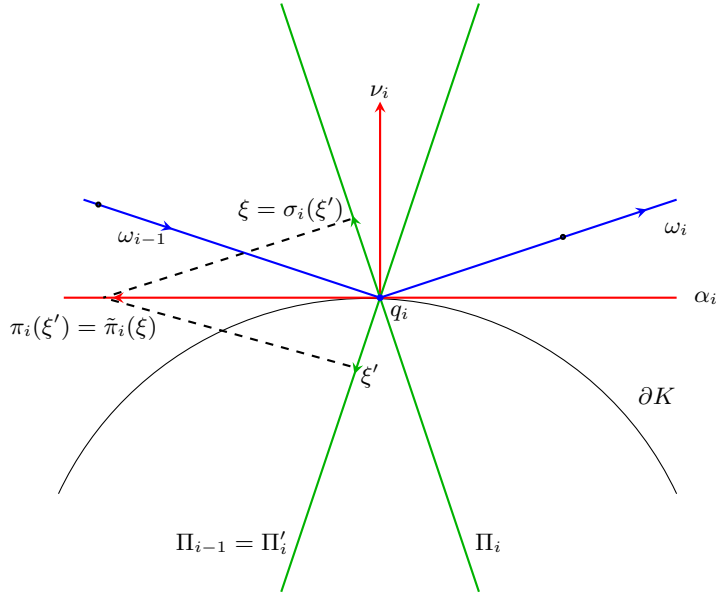
Figure 2: The projections $\pi_i$ and $\tpi_i$
\end{center}

\bs

We will now write down some of the above formulae a bit differently.
Let 
$$\tpi_i : \Pi_i \longrightarrow \alpha_i$$
be the {\it projection} along the vector $\omega_{i}$. For the symmetry $\sigma_i: \R^3 \longrightarrow \R^3$ through the plane
$\alpha_i = T_{q_i}(\dk)$ we have
$$\sigma_i(\xi) = \xi - 2 \langle \xi, \nu_i \rangle \, \nu_i ,$$
for all $\xi \in \R^3$.
Since $\theta_i = 2 \langle \omega_i , \nu_i\rangle = - 2 \langle \omega_{i-1} , \nu_i\rangle$, we have
\begin{equation}
\la \tpsi_i (\xi) , \eta \ra =   \theta_i \,\la G_i(\tpi_i(\xi)) , \tpi_i(\eta) \ra
\end{equation}
for all $\xi, \eta \in \Pi_i$, where we used (3.4) and the fact that 
$\pi_i(\xi') = \pi_i (\sigma_i(\xi')) = \tpi_i(\xi)$ for $\xi' \in \Pi_i'$ and $\xi = \sigma_i(\xi') \in \Pi_i$.
The latter follows from the formula for $\sigma_i$ and 
\begin{equation}
\pi_i(\xi') = \xi' - \frac{\langle \xi', \nu_i\rangle}{\langle \omega_{i-1}, \nu_i\rangle}  \, \omega_{i-1} 
\quad, \quad \tpi_i(\xi) = \xi - \frac{\langle \xi, \nu_i\rangle}{\langle \omega_{i}, \nu_i\rangle}  \, \omega_{i} 
\end{equation}
for all $\xi' \in \Pi_i'$ and all $\xi \in \Pi_i$ (or just use simple geometry).

We now continue with the proof of Theorem 1.1.

Assume that the centres $O_i$ of the balls $K_i$ are all in the plane $\Gamma$. It is easy to see that under this assumption 
all billiard trajectories in the non-wandering set $M_0$ are lying in $\Gamma$. Thus, 
the whole billiard trajectory  $q_0, q_1, \ldots$ lies in the plane $\Gamma$, so {\bf all vectors $\omega_i$ and
normals $\nu_i$ are in this plane, so perpendicular to $n = (0,0,1)$}. All planes $\Pi_i$ are then 
perpendicular to $\Gamma$, so from the identifications we made, we can assume $n \in \Pi_i$ for all $i$.
(Identifying all $\Pi_i$ with the corresponding parallel 2-dimensional planes in $\R^3$ through $0$.)
All tangent planes $\alpha_i$ are also perpendicular to $\Gamma$ and each of them is identified with a corresponding
parallel plane in $\R^3$ through $0$.

Now the second fundamental form $G_i$ acts in a very simple way on $\alpha_i$, namely
$$G_i(\xi) = \kappa_i \, \xi $$ 
for all $\xi \in \alpha_i$, where $\kappa_i = 1/r_i$. Now (3.8) implies
\begin{equation}
\la \tpsi_i(\xi), \eta\ra = \kappa_i\, \theta_i \, \la \tpi_i(\xi), \tpi_i(\eta) \ra
\end{equation}
for all $\xi, \eta \in \Pi_i$. For each $i$ denote by $E_i$ {\it the $1$-dimensional  subspace of $\Pi_i$ which is perpendicular}
to $\Gamma$. This is just a line through $0$ in $\Pi_i$ which contains the normal vector $n = (0,0,1)$ to $\Gamma$.
Assuming $\Pi_i \neq \alpha_i$ (which we do), it follows that 
$$E_i = \Pi_i \cap \alpha_i .$$
Since $E_i \perp \Gamma$, it is perpendicular to all vectors in the plane. In particular,
$E_i \perp \nu_i$, so $\tpi_i(\xi) = \xi$ for all $\xi \in E_i$. Now (3.10) yields
$\la \tpsi_i(\xi), \eta\ra =  \kappa_i\, \theta_i \, \la \xi, \eta \ra$
for all $\xi, \eta \in E_i$.  Therefore the restriction of $\tpsi_i$ to the subspace $E_i$ is
\begin{equation}
(\tpsi_i)_{|E_i} = \kappa_i\, \theta_i \, I ,
\end{equation}
where $I$ is the identity operator. A similar formula holds for $\psi_i$, so now for the operator
$$\Pm_\gamma : E_0 \times E_0 \longrightarrow E_m \times E_m$$
we obtain the following matrix representation
\be 
(\Pm_\gamma)_{| E_0 \times E_0} =
 \left(
\begin{array}{cc}
s_m & 0\\
0  & s_m
\end{array}
\right)
\left(
\begin{array}{cc}
I & d_m I\\
\kappa_m \,\theta_m \, I & (1 + \kappa_m\, d_m \theta_m) I
\end{array}
\right)
\cdots
\left(
\begin{array}{cc}
I & d_1 I\\
\kappa_1 \,\theta_1 \, I  & (1 + \kappa_1 \, d_1 \theta_1) I
\end{array}
\right) .
\ee

Set $F_i = \Pi_i \cap \Gamma$. This is then an one-dimensional subspace of $\Pi_i$ contained in $\Gamma$. 
Since $\sigma_i(\Pi_{i-1}) = \Pi_i$, we have $\sigma_i(F_{i-1}) = F_i$, and so $s_m(F_0) = F_m$.

Clearly,
$\omega_i \perp F_i$, since $\omega_i \perp \Pi_i$. Let $p_i$ be the projection of $\nu_i$ in $\Pi_i$, i.e.
$$p_i = \nu_i - \langle \nu_i, \omega_i\rangle\, \omega_i .$$
Then $p_i \in {\mbox{\rm span}} (\nu_i, \omega_i) \subset \Gamma$ and moreover $p_i \in F_i$. Assuming that the billiard trajectory
generated by $\tx_0$ is never perpendicular to $\dk$ (by the choice of $\tx_0$), i.e. $\nu_i \neq \omega_i$ for all $i$, 
we have $p_i \neq 0$. Set
$$f_i = \frac{p_i}{\|p_i\|} .$$
Since $f_i \in F_i$ and $F_i$ is one-dimensional, $F_i = \{ s \, f_i : s \in \R\}$.

We will now calculate $\tpsi_i(f_i)$ using (3.10). First, (3.9) implies
\begin{eqnarray*}
\tpi_i(p_i) 
& = & p_i - \frac{\langle p_i , \nu_i\rangle}{\langle \omega_i , \nu_i\rangle} \, \omega_i\\
& = & \nu_i - \langle \nu_i, \omega_i\rangle\, \omega_i 
- \frac{\langle  \nu_i - \langle \nu_i, \omega_i\rangle\, \omega_i , \nu_i\rangle}{\langle \omega_i , \nu_i\rangle} \, \omega_i\\
& = & \nu_i - \langle \nu_i, \omega_i\rangle\, \omega_i - \frac{1}{\langle \omega_i , \nu_i\rangle}\, \omega_i
+ \langle \nu_i, \omega_i\rangle\, \omega_i \\
& = & \nu_i - \frac{1}{\langle \omega_i , \nu_i\rangle}\, \omega_i .
\end{eqnarray*}
Thus, 
$$\|\tpi_i(p_i)\|^2 = 1-2 + \frac{1}{\langle \omega_i, \nu_i\rangle^2} 
 = \frac{1 - \langle \omega_i , \nu_i\rangle^2}{\langle \omega_i, \nu_i\rangle^2} .$$
Now (3.10) yields
$$\langle \tpsi_i (p_i) , p_i \rangle 
= \kappa_i\, \theta_i\, \|\tpi_i(p_i)\|^2 = \kappa_i \, \theta_i \, \frac{1 - \langle \omega_i , \nu_i\rangle^2}{\langle \omega_i, \nu_i\rangle^2}.$$
Since $\|p_i\|^2 = 1 - \langle \omega_i , \nu_i\rangle^2$, it follows that
\begin{equation}
\langle \tpsi_i ( f_i) , f_i \rangle = \frac{\kappa_i\, \theta_i}{\langle \omega_i, \nu_i\rangle^2} .
\end{equation}

We will now show that $\tpsi_i (f_i) \in F_i$. For any $\xi \in E_i$, using (3.10) again and the fact that $\tpi_i(\xi) = \xi \perp \nu_i$ we get
\begin{eqnarray*}
\langle \tpsi_i (f_i) , \xi \rangle
& = & \kappa _i\, \theta_i \, \langle \tpi_i(f_i) , \xi\rangle\\
& = & \kappa_i \, \theta_i \, \left\langle  f_i - \frac{\langle f_i, \nu_i\rangle}{\langle \omega_i , \nu_i\rangle}\, \omega_i , \xi  \right\rangle
 =   \kappa_i \, \theta_i \, \langle f_i, \xi \rangle .
\end{eqnarray*}
This yields $\langle \tpsi_i (f_i) - \kappa_i \, \theta_i \,  f_i , \xi \rangle = 0$ for all $\xi \in E_i$, so 
$\tpsi_i (f_i) - \kappa_i \, \theta_i \,  f_i \in E_i^\perp$, i.e. it is in the plane $\Gamma$. 
Since, $\tpsi_i (f_i) - \kappa_i \, \theta_i \,  f_i \in \Pi_i$,
we must have $\tpsi_i (f_i) - \kappa_i \, \theta_i \,  f_i \in F_i$, i.e. $\tpsi_i (f_i) - \kappa_i \, \theta_i \,  f_i $ is parallel to $f_i$.
In particular, $\tpsi_i (f_i) = t\, f_i$ for some $t \in \R$. Now (3.13) implies $t = \frac{\kappa_i\, \theta_i}{\langle \omega_i, \nu_i\rangle^2} $.
Hence
\begin{equation}
\tpsi_i (f_i)  = \kappa_i\, \ttheta_i \,  f_i ,
\end{equation}
where
$$\ttheta_i = L_i \theta_i \geq \theta_i \quad, \quad L_i = \frac{1}{\langle \omega_i, \nu_i\rangle^2} \geq 1$$
for all $i$.

It follows from the Assumptions in Sect. 1.1 that  there exists a global constant 
$c_2 < 1$ such that  if the reflection points $q_{i-1}$ and $q_{i+1}$  (i.e. the ones before and after $q_i$) belong to
different components of $K$, then $\langle \omega_i, \nu_i\rangle \leq c_2$. Therefore in such a case we have
\begin{equation}
L_i \geq L = \frac{1}{c_2^2} > 1 ,
\end{equation}
where $L > 1$ defined here  is a global constant.

\bs

Now we get a similar formula to (3.12)  for the operator
$$\Pm_\gamma : F_0 \times F_0 \longrightarrow F_m \times F_m ,$$
namely we have the following matrix representation
\be 
(\Pm_\gamma)_{| F_0 \times F_0} =
 \left(
\begin{array}{cc}
s_m & 0\\
0  & s_m
\end{array}
\right) \cdot
\left(
\begin{array}{cc}
I & d_m I\\
\kappa_m \,\ttheta_m \, I & (1 + \kappa_m \, d_m \ttheta_m) I
\end{array}
\right)
\cdots
\left(
\begin{array}{cc}
I & d_1 I\\
\kappa_1 \,\ttheta_1 \, I  & (1 + \kappa_1 \, d_1 \ttheta_1) I
\end{array}
\right) .
\ee

\ms

Set for brevity
$$a_i = \kappa_i\, \theta_i \quad, \quad \ta_i = \kappa_i\, \ttheta_i  \geq a_i\quad, \quad
b_i = 1 + \kappa_i\, d_i\, \theta_i \quad, \quad \tb_i = 1 + \kappa_i\, d_i\, \ttheta_i \geq b_i.$$

\ms

\noindent
{\bf Lemma 3.2.} {\it For every $j = 1,2, \ldots, m$ there exist functions
$$ \Deone_j = \Deone_j (\{\kappa_i\}_{i=1}^j; \{d_i\}_{i=1}^j ; \{\theta_i\}_{i=1}^j) \quad, \quad
 \Detwo_j = \Detwo_j (\{\kappa_i\}_{i=1}^j; \{d_i\}_{i=1}^j ; \{\theta_i\}_{i=1}^j) ,$$
$$ \Dethree_j = \Dethree_j (\{\kappa_i\}_{i=1}^j; \{d_i\}_{i=1}^j ; \{\theta_i\}_{i=1}^j)
\quad, \quad  \Defour_j = \Defour_j (\{\kappa_i\}_{i=1}^j; \{d_i\}_{i=1}^j ; \{\theta_i\}_{i=1}^j)$$
which are polynomials of $\kappa_1, \ldots, \kappa_j; d_1, \ldots, d_j ; \theta_1, \ldots, \theta_j$ with positive coefficients
in front of all terms such that
\be 
\left(
\begin{array}{cc}
I & d_j I\\
a_j \, I & b_j\, I
\end{array}
\right)
\cdots
\left(
\begin{array}{cc}
I & d_1 I\\
a_1 \, I  & b_1\, I
\end{array}
\right) 
= \left(
\begin{array}{cc}
\Deone_j\, I & \Detwo_j\,  I\\
\Dethree_j \, I  & \Defour_j\, I 
\end{array}
\right) .
\ee
}

\medskip

\noindent
{\it Proof.} For $j = 1$ the statement is obviously true. Assuming the statement is true for some $j\geq 1$, we get
\begin{eqnarray*}
&      &\left(
\begin{array}{cc}
I & d_{j+1} I\\
a_{j+1} \, I & b_{j+1}\,  I
\end{array}
\right)\;
\left(
\begin{array}{cc}
I & d_j I\\
a_j \, I & b_j\, I
\end{array}
\right)
\cdots \left(
\begin{array}{cc}
I & d_1 I\\
a_1  \, I  & b_1\, I
\end{array}
\right) \\
& = & 
\left(
\begin{array}{cc}
I & d_{j+1} I\\
a_{j+1}  \, I & b_{j+1}\, I
\end{array}
\right)\;
\left(
\begin{array}{cc}
\Deone_j\, I & \Detwo_j\,  I\\
\Dethree_j \, I  & \Defour_j\, I 
\end{array}
\right) \\
& = & \left(
\begin{array}{cc}
(\Deone_j + d_{j+1}\, \Dethree_j)\, I & (\Detwo_j + d_{j+1}\, \Defour_j) \,  I\\
(a_{j+1}\, \Deone_j + b_{j+1}\,\Dethree_j )\, I  
& (a_{j+1} \Detwo_j  + b_{j+1} \,\Defour_j)\, I 
\end{array}
\right) 
 =  \left(
\begin{array}{cc}
\Deone_{j+1}\, I & \Detwo_{j+1}\,  I\\
\Dethree_{j+1} \, I  & \Defour_{j+1}\, I 
\end{array}
\right) ,
\end{eqnarray*}
where we set
$$\Deone_{j+1} = \Deone_j + d_{j+1}\, \Dethree_j \quad, \quad  \Detwo_{j+1} = \Detwo_j + d_{j+1}\, \Defour_j  ,$$
$$\Dethree_{j+1} = a_{j+1}\, \Deone_j + b_{j+1} \,\Dethree_j \quad, \quad 
\Defour_{j+1} =  a_{j+1} \Detwo_j + b_{j+1}\,\Defour_j .$$
This proves the lemma.
\endofproof

\bs

Setting
$$ \tDeone_j = \Deone_j (\{\kappa_i\}_{i=1}^j; \{d_i\}_{i=1}^j ; \{\ttheta_i\}_{i=1}^j) \quad, \quad
 \tDetwo_j = \Detwo_j (\{\kappa_i\}_{i=1}^j; \{d_i\}_{i=1}^j ; \{\ttheta_i\}_{i=1}^j) ,$$
$$ \tDethree_j = \Dethree_j (\{\kappa_i\}_{i=1}^j; \{d_i\}_{i=1}^j ; \{\ttheta_i\}_{i=1}^j)
\quad, \quad  \tDefour_j = \Defour_j (\{\kappa_i\}_{i=1}^j; \{d_i\}_{i=1}^j ; \{\ttheta_i\}_{i=1}^j) ,$$
as in the proof of Lemma 3.2, we get the following
\be 
\left(
\begin{array}{cc}
I & d_j I\\
\ta_j \, I & \tb_j\, I
\end{array}
\right)
\cdots
\left(
\begin{array}{cc}
I & d_1 I\\
\ta_1 \, I  & \tb_1\, I
\end{array}
\right) 
= \left(
\begin{array}{cc}
\tDeone_j\, I & \tDetwo_j\,  I\\
\tDethree_j \, I  & \tDefour_j\, I 
\end{array}
\right) .
\ee

Now combining (3.2) with the representations (3.12) and (3.17), and also (3.16) and (3.18) we obtain
the following consequence concerning Lyapunov exponents.

\bs

\noindent
{\bf Corollary 3.3.} {\it For all $v \in E_0$ and $\tv \in F_0$ with $\|v\| = \|\tv\| = 1$ we have
\be
\lambda(v) = \lim_{m\to\infty} \frac{1}{m} \log  \Defour_m  ,
\ee
and}
\be
\lambda(\tv) = \lim_{m\to\infty} \frac{1}{m} \log  \tDefour_m  .
\ee

\ms

\noindent
{\it Proof.} For $v \in E_0$, (3.12) and (3.17) imply
\begin{equation}
\Pm_\gamma\, \left(
\begin{array}{c}
0\\
v
\end{array}
\right) =
\left(
\begin{array}{cc}
s_m & 0\\
0 & s_m
\end{array}
\right) \cdot 
\left(
\begin{array}{cc}
\Deone_m\, I & \Detwo_m\,  I\\
\Dethree_m \, I  & \Defour_m\, I 
\end{array}
\right) \cdot 
\left(
\begin{array}{c}
0\\
v
\end{array}
\right) =
\left(
\begin{array}{c}
\Detwo_m\, s_m(v)\\
\Defour_m \, s_m(v)
\end{array}
\right) .
\end{equation}
Therefore
$$\pr_2 \left(
\Pm_\gamma\, \left(
\begin{array}{c}
0\\
s_m(v)
\end{array}
\right) 
\right) = \Defour_m \, s_m(v) .$$
Since $\|s_m(v)\| = \|v\|$, now (3.19) follows from the above and  (3.2).

In a similar way one proves (3.20).
\endofproof

\bs


Recall from the proof of Lemma 3.2 that
\begin{equation}
\Detwo_{j+1} = \Detwo_j + d_{j+1}\, \Defour_j \quad, \quad 
\Defour_{j+1} =  a_{j+1}\,  \Detwo_j + b_{j+1}\, \,\Defour_j .
\end{equation}
Similarly,
\begin{equation}
\tDetwo_{j+1} = \tDetwo_j + d_{j+1}\, \tDefour_j \quad, \quad 
\tDefour_{j+1} =  \ta_{j+1}\,  \Detwo_j + \tb_{j+1}\, \,\Defour_j .
\end{equation}
We set $\Detwo_0 = \tDetwo_0 = \Defour_0 = \tDefour_0 = 1$.

\bs

\noindent
{\bf Lemma 3.4.} {\it For all $j \geq 1$ we have}
\begin{equation}
\Detwo_j \leq \Defour_j \quad, \quad \tDetwo_j \leq \tDefour_j .
\end{equation}

\ms

\noindent
{\it Proof.} We will prove the first half of the above; the other is similar.

For $j = 1$ the statement is true. Assume it is true for some $j \geq 1$.
Then (3.22) implies
$$\Detwo_{j+1} = \Detwo_j + d_{j+1}\, \Defour_j \leq 2 \Defour_j ,$$
since $d_{j+1} \leq 1$. This and $b_{j+1} = 1 + \kappa_{j+1} \theta_{j+1} d_{j+1} \geq 2$ (see the Assumptions in Sect. 1) imply
$$\Defour_{j+1} =  a_{j+1}\,  \Detwo_j + b_{j+1}\, \,\Defour_j  \geq b_{j+1}\, \,\Defour_j \geq
2 \Defour_j \geq \Detwo_{j+1} .$$
By induction $\Detwo_i \leq \Defour_i$ for all $i \geq 1$.
\endofproof

\bs

Here we used $\kappa_i = 1/r_i \geq 1/r_0$, so
$\kappa_i \theta_i \geq c_2/r_0 > 2$ for all $i$ (see the Assumptions in Sect. 1).
From the definitions, $\tb_i \geq b_i$ for all $i$. 
Since $L_i \geq L > 1$ whenever $q_{i-1}, q_i, q_{i+1}$ belong to distinct balls, it follows that
$$\tb_i = 1 + \kappa_i d_i \ttheta_i \geq  1 + \kappa_i d_i \theta_i L$$
for such $i$. Set 
$$I = \{ i : 1 \leq i \leq m\:, \tb_i > b_i \} .$$
Then
$$\tb_i - b_i \geq \kappa_i\, d_i \theta_i (L-1) \geq c = \con > 0$$
for $ \in I$. We choose $c \leq 1$ so that then $\tb_i \geq b_i \geq c$ for all $i$.

Similarly, 
$$\ta_i - a_i =  \kappa_i\, \theta_i  (L-1) \geq c = \con > 0 ,$$
assuming the global constants $\kappa_i = 1/r_i$ are sufficiently large (see the Assumptions in Sect. 1).
Since $a_i$ and $\ta_i$ are uniformly bounded above and below by global positive constants,
there exists a global constant $C > 0$ so that 
$$\ta_i  \geq a_i + c \geq a_i (1+ c/C).$$
Choose a small global constant $s > 0$, with $1+ c/C \geq e^s$. Then $\ta_i \geq e^s\, a_i$ for all $i$. In a similar way
we deal with $\tb_i$ and $b_i$, so we get
\begin{equation}
\ta_i \geq e^s\, a_i \quad, \quad \tb_i \geq e^s \, b_i
\end{equation}
for all $i \in I$. For $i \notin I$ we simply have $\ta_i = a_i$ and $\tb_i = b_i$.
We assume $s > 0$ is sufficiently small so that $e^s \leq 1+ d_0$.

It is clear from (3.22) and (3.23) that $\Detwo_j$ and $\tDetwo_j$ grow exponentially with $j$.
E.g. (3.22) and (3.24) give
$$\Detwo_{j+1} = \Detwo_j + d_{j+1} \Defour_j \geq (1+ d_{j+1}) \Detwo_j $$
for all $j$, so 
$\Detwo_j \geq \prod_{i=1}^j (1 + d_i) \geq (1+d_0)^j.$
In a similar way one can see that $\Defour_j$ and $\tDefour_j$ also grow exponentially with $j$. However we need more precise estimates,
and we need to compare the growth of $\tDefour_j$ with that of $\Defour_j$.

Denote by $i_0 > 0$ the largest integer so that $i_0 \notin I$. Then $i_0 < k$. Set 
\begin{equation}
u_i = v_i = 1\quad , \quad 0 \leq i \leq i_0 .
\end{equation}

\bs

\noindent
{\bf Lemma 3.5.} {\it There exist sequences $\{t_i\}_{i\geq i_0}^\infty$ and $\{s_i\}_{i\geq i_0}^\infty$ of non-negative numbers so that
\begin{equation}
u_{i+1} = e^{t_{i+1}} \, u_i \quad, \quad v_{i} = e^{s_i} \, u_i 
\end{equation}
for all $i \geq i_0$, $u_i \geq 1$ for all $i \geq 1$, and
\begin{equation}
\tDetwo_{i} \geq u_i\, \Detwo_i
\end{equation}
and
\begin{equation}
\tDefour_{i} \geq v_i\, \Defour_i
\end{equation}
for all $i \geq i_0$. Moreover we can choose these sequences so that 
\begin{equation}
0 \leq t_{i+1} < s_i
\end{equation}
for $i \geq i_0$ and }
\begin{equation}
s_i \geq s/2 \quad, \quad i \in I .
\end{equation}

\ms

\noindent
{\it Proof.} We will do the construction by induction. For $0 \leq i \leq i_0$ we use (3.26).

Let us check (3.28) and (3.29) for  $i  = i_0$. By (3.22), (3.23)  and $u_{i_0} = v_{i_0} = 1$ we get
\begin{eqnarray*}
\tDetwo_{i_0}
& =      & \tDetwo_{i_0-1} + d_{i_0}\, \tDefour_{i_0-1}
 \geq  \Detwo_{i_0-1} + d_{i_0}\,  \Defour_{i_0-1} = \Detwo_{i_0} = u_{i_0}\, \Detwo_{i_0},
\end{eqnarray*}
and
\begin{eqnarray*}
\tDefour_{i_0}
& =      & \ta_{i_0}\, \tDetwo_{i_0-1} + \tb_{i_0}\, \tDefour_{i_0-1}
=  a_{i_0}\, \tDetwo_{i_0-1} + b_{i_0}\, \tDefour_{i_0-1} \geq   a_{i_0}\, \Detwo_{i_0-1} + b_{i_0}\, \Defour_{i_0-1} = v_{i_0} \, \Defour_{i_0} .
\end{eqnarray*}

Next, recall that for $i = i_0$ we have $i_0 + 1 \in I$. Set $t_{i_0+1} = 0$ and $u_{i_0 +1} = 1 = e^{t_{i_0+1}}\, u_{i_0}$. Then clearly
$$\tDetwo_{i_0 +1} \geq u_{i_0 + 1}\, \Detwo_{i_0 + 1} ,$$
so (3.28) holds with $i$ replaced by $i_0 + 1$. 
Now set $s_{i_0 + 1} = s$ and $v_{i_0 + 1} = e^{s_{i_0 + 1}}u_{i_0 + 1} = e^{s}$. Then (3.23), (3.22), $i_0 + 1 \in I$ and (3.25) imply
\begin{eqnarray*}
\tDefour_{i_0 + 1}
& =      & \ta_{i_0 + 1}\, \tDetwo_{i_0} + \tb_{i_0+1}\, \tDefour_{i_0}
 \geq  e^s\, a_{i_0 + 1} \,\Detwo_{i_0} + e^s\, b_{i_0 + 1}\,  \Defour_{i_0}  = e^s \, \Defour_{i_0 + 1} ,
\end{eqnarray*}
so (3.29) holds with $i$ replaced by $i_0 + 1$. 

Let us now show that we can choose $t_{i_0+2} > 0$ so that setting $u_{i_0 + 2} = e^{t_{i_0 + 2}}\, u_{i_0+1} = e^{t_{i_0 + 2}}$,
(3.28) holds with $i$ replaced by $i_0 + 2$. First, using the above and (3.23) we get
$$\tDetwo_{i_0 + 2} = \tDetwo_{i_0 + 1} + d_{i_0 + 2}\, \tDefour_{i_0 + 1}
\geq \Detwo_{i_0 + 1} + d_{i_0 + 2} e^s \, \Defour_{i_0+1} .$$
Now, using (3.22), we have $d_{i_0+2}\, \Defour_{i_0 + 1} = \Detwo_{i_0 + 2} - \Detwo_{i_0 + 1}$, so the above implies
\begin{equation}
\tDetwo_{i_0 + 2} \geq \Detwo_{i_0 + 1} + e^s \, (\Detwo_{i_0 + 2} - \Detwo_{i_0 + 1})
= e^s \, \Detwo_{i_0 + 2} - (e^s - 1) \Detwo_{i_0 + 1} .
\end{equation}
However, by (3.22) and (3.24),
$$\Detwo_{i_0 + 2} = \Detwo_{i_0 + 1} + d_{i_0 + 2}\, \Defour_{i_0 + 1} \geq (1+ d_0)\, \Detwo_{i_0 + 1} ,$$
therefore
$\Detwo_{i_0 + 1} \leq \frac{\Detwo_{i_0 + 2}}{1+d_0} .$
Combining this with (3.32) implies
$$\tDetwo_{i_0 + 2} \geq \left(e^s - \frac{{e^s}-1}{1+d_0} \right)\, \Detwo_{i_0 + 2} 
=  \left(\frac{d_0\, e^s + 1}{1+d_0} \right)\, \Detwo_{i_0 + 2} .$$
Since $\frac{d_0\, e^s + 1}{1+d_0}  > 1$, there exists $t_{i_0 + 2} > 0$ so that 
$$0 < t_{i_0 + 2} < s = s_{i_0+1} \quad, \quad \frac{d_0 {e^s} + 1}{1+d_0} \geq e^{t_{i_0 + 2}} .$$
Choose an arbitrary $t_{i_0+2} > 0$ with the above and set
$u_{i_0 + 2} = e^{t_{i_0 + 2}} = e^{t_{i_0 + 2}} \, u_{i_0 + 1}$. Then (3.28) holds with $i$ replaced by $i_0 + 2$.

Assume that $u_i$, $t_i$, $v_i$, $s_i$ have been constructed for some $i \geq i_0+1$ so that (3.27)-- (3.30) hold for $i$.
We will now define $u_{i+1}$, $t_{i+1}$, $v_{i+1}$, $s_{i+1}$. For the first two we will essentially repeat the previous argument.

Using (3.23) and (3.27), (3.29) for $i$, we get
$$\tDetwo_{i + 1} = \tDetwo_{i } + d_{i + 1}\, \tDefour_{i}
\geq u_i\, \Detwo_{i} + d_{i + 1} \, v_i\, \Defour_{i} = u_i\, \Detwo_{i} + d_{i + 1} e^{s_i} u_i\, \Defour_{i} .$$
Now, using (3.22), we have $d_{i + 1}\, \Defour_{i} = \Detwo_{i + 1} - \Detwo_{i}$, so the above yields
\begin{equation}
\tDetwo_{i + 1} \geq u_i \, \left(\Detwo_{i} + e^{s_i} \, (\Detwo_{i + 1} - \Detwo_{i}) \right)
= u_i\, \left( e^{s_i} \, \Detwo_{i + 1} - (e^{s_i} - 1) \Detwo_{i} \right).
\end{equation}
From (3.22) and (3.24),
$\Detwo_{i + 1} = \Detwo_{i} + d_{i + 1}\, \Defour_{i} \geq (1+ d_0)\, \Detwo_{i} ,$
so
$\Detwo_{i} \leq \frac{\Detwo_{i + 1}}{1+d_0} .$
Using this in (3.33) gives
$$\tDetwo_{i + 1} \geq u_i\, \left(e^{s_i} - \frac{{e^{s_i}}-1}{1+d_0} \right)\, \Detwo_{i + 1} 
=  u_i\, \left(\frac{d_0\, e^{s_i} + 1}{1+d_0} \right)\, \Detwo_{i + 1} .$$
Since $\frac{d_0\, e^{s_i} + 1}{1+d_0}  > 1$, there exists $t_{i+1} > 0$ so that
\begin{equation}
0 < t_{i + 1} < s_i \quad, \quad \frac{d_0 {e^{s_i}} +1}{1+d_0} \geq e^{t_{i + 1}} .
\end{equation}
Choose $t_{i+1} > 0$ with these properties and set 
$u_{i + 1} = e^{t_{i + 1}}\, u_i$. Then  the above gives $\tDetwo_{i + 1} \geq u_{i+1}\, \Detwo_{i+1}$, so
 (3.28) holds with $i$ replaced by $i+1$. We will impose another condition on $t_{i+1}$ later.

Next, we will define $s_{i+1}$ and $v_{i+1} = e^{s_{i+1}} \, u_{i+1}$ so that (3.29) holds with $i$ replaced by $i+1$. 
There are two cases to consider.

First assume that $i+1 \notin I$. Then (3.23), (3.27) and the inductive assumptions yield
\begin{eqnarray}
\tDefour_{i+1}
& =      & \ta_{i+1}\, \tDetwo_i + \tb_{i+1}\, \tDefour_i  =  a_{i+1}\, \tDetwo_i + b_{i+1}\, \tDefour_i\nonumber\\
& \geq & a_{i+1}\, u_i\, \Detwo_i + b_{i+1}\, v_i\, \Defour_i = a_{i+1}\, u_i\, \Detwo_i + b_{i+1}\, e^{s_i}\, u_i \, \Defour_i\nonumber\\
& =      & e^{t_{i+1}}\, u_i \, \left(a_{i+1}\, e^{-t_{i+1}}\,  \Detwo_i + b_{i+1}\, e^{s_i - t_{i+1}}\, \Defour_i \right) \nonumber\\
& =      & u_{i+1} \, \left(a_{i+1}\, e^{-t_{i+1}}\,  \Detwo_i + b_{i+1}\, e^{s_i - t_{i+1}}\, \Defour_i \right) .
\end{eqnarray}
Thus, we need to define $v_{i+1} = e^{s_{i+1}} \, u_{i+1}$ so that
$$u_{i+1} \, \left(a_{i+1}\, e^{-t_{i+1}}\,  \Detwo_i + b_{i+1}\, e^{s_i - t_{i+1}}\, \Defour_i \right) 
\geq v_{i+1} \Defour_{i+1} = e^{s_{i+1}} \, u_{i+1} \, \left(a_{i+1}\,  \Detwo_i + b_{i+1}\, \Defour_i \right) . $$
The latter is equivalent to
\begin{equation}
b_{i+1}\, \Defour_i \, \left(e^{s_{i} - t_{i+1}} - e^{s_{i+1}} \right) \geq a_{i+1}\, \Detwo_i \, \left( e^{s_{i+1}} - e^{-t_{i+1}}\right) .
\end{equation}
To simplify the above, notice that
$$\frac{b_{i+1}}{a_{i+1}} = \frac{1 + d_{i+1} \, a_{i+1}}{a_{i+1}}
= \frac{1}{a_{i+1}} + d_{i+1} \geq d_0 . $$
This and (3.24) imply
$$\frac{b_{i+1}}{a_{i+1}}\cdot \frac{\Defour_i}{\Detwo_i} \geq d_0  .$$
Hence to satisfy (3.36) it is enough to have
$$ d_0 \,\left(\, e^{s_{i} - t_{i+1}} - e^{s_{i+1}}\right) \geq e^{s_{i+1}} - e^{-t_{i+1}} ,$$
that is
\begin{equation}
d_0 \, \left(e^{s_{i} - t_{i+1} - s_{i+1}} - 1 \right) \geq 1 - e^{-(t_{i+1} + s_{i+1})} .
\end{equation}
Clearly $1 - e^{-(t_{i+1} + s_{i+1})}$ is close to $0$ when $t_{i+1} + s_{i+1}$ is close to $0$, while $d_0 \, \left(e^{s_{i} - t_{i+1} - s_{i+1}} - 1 \right)$
is close to $d_0 (e^{s_i} -1) > 0$. Thus, (3.37) holds when $t_{i+1} > 0$ and $s_{i+1} > 0$ are chosen sufficiently close to $0$.
 So, to finish with the case under consideration we just choose arbitrary $t_{i+1} > 0$ and $s_{i+1} > 0$ with 
$0 < t_{i+1} + s_{i+1} <  s_i$ so that (3.34) and (3.37) hold. Then (3.36) holds as well and as observed earlier, it
implies $\tDefour_{i+1} \geq v_{i+1} \, \Defour_{i+1}$, i.e. (3.29) holds with $i$ replaced by $i+1$.

Next, consider the case $i+1 \in I$. Then (3.23) and the inductive assumptions give
\begin{eqnarray}
\tDefour_{i+1}
& =      & \ta_{i+1}\, \tDetwo_i + \tb_{i+1}\, \tDefour_i  \geq e^s\, \left(  a_{i+1}\, \tDetwo_i + b_{i+1}\, \tDefour_i\right)\nonumber\\
& \geq & e^s\, \left(a_{i+1}\, u_i\, \Detwo_i + b_{i+1}\, v_i\, \Defour_i \right)
= e^s\, \left( a_{i+1}\, u_i\, \Detwo_i + b_{i+1}\, e^{s_i}\, u_i \, \Defour_i \right) \nonumber\\
& \geq & e^{s}\, u_i \, \left(a_{i+1}\,  \Detwo_i + b_{i+1}\, \Defour_i \right) 
= e^s\, u_i\, \Defour_{i+1} \nonumber\\
& =      & e^{s- t_{i+1}}\, u_{i+1}\, \Defour_{i+1} .
\end{eqnarray}
In this case we just take $t_{i+1} < s/2$ so that (3.34) holds and set $s_{i+1} = s - t_{i+1} > s/2$ and $v_{i+1} = e^{s_{i+1}}\, u_{i+1}$.
Then (3.29), (3.30) and (3.31) hold with $i$ replaced by $i+1$.

This completes the induction.
\endofproof

\bs

To finish with the proof of Theorem 1.1 we will use the fact that $z \in \Sigma_0$, the set defined by (2.6), and the number $k$
with (2.4); see Sect. 2.3 above.
We can assume that $m = k m'$ for some $m' \geq 6$. Then, by Corollary 2.3, for $\ell = |I|$ we have $\ell \geq m/2 - k \geq m/3$.
For $v \in E_0$ and $\tv \in F_0$ with $\|v\| = \|\tv\| = 1$,
(3.19), (3.20),  (3.29) and  (3.31) imply
\begin{eqnarray*}
\lambda(\tv) 
& =      & \lim_{m\to\infty} \frac{1}{m} \log  \tDefour_m  \geq
\lim_{m\to\infty} \frac{1}{m} \log  (e^{\ell \, s/2}\, \Defour_m  )\\
& =      & \lim_{m\to\infty} \frac{1}{m} \left(\frac{\ell \,s}{2} + \log \Defour_m  \right)\\
& \geq & \frac{s}{6} + \lim_{m\to\infty} \frac{1}{m} \log  \Defour_m  = \frac{s}{6} + \lambda(v) .
\end{eqnarray*}

The above shows that we must have $\lambda_1 = \lambda(\tv)$ for $\tv \in F_0$ and
$\lambda_2 = \lambda(v)$ for $v \in E_0$ and also that $\lambda_1 > \lambda_2$.

This completes the proof of Theorem 1.1.
\endofproof

\bs

{\footnotesize

\end{document}